\documentclass[12pt]{article}
\usepackage{latexsym, amsfonts}
\usepackage{epsfig}

\newtheorem{thm}{Theorem}[section]
\newtheorem{prop}{Proposition}[section]

\renewcommand{\frac}[2]{\textstyle {#1\over #2}\displaystyle}

\def\ex{\mbox{\rm E\,}}
\def\var{\mbox{\rm Var\,}}

\def\I{\mbox{\bf 1}}

\begin{document}

\title{Estimating  the Structural Distribution\\
Function of Cell Probabilities}

\author{Bert   van Es, Chris A.J. Klaassen \\
{\normalsize Korteweg-de Vries Institute for Mathematics}\\
{\normalsize University of Amsterdam}\\
\and
 Robert M. Mnatsakanov
 \thanks{Financed by INTAS-97-Georgia-1828} \\
{\normalsize West Virginia University, Morgantown}\\
{\normalsize A. Razmazde Mathematical Institute, Tbilisi}
}

\maketitle

\begin{abstract}
We consider estimation of the structural distribution function of
the cell probabilities of a multinomial sample in situations where
the number of cells is large. We review the performance of the
natural estimator, an estimator based on grouping the cells and a
kernel type estimator. Inconsistency of the natural estimator and
weak consistency of the other two estimators is derived by
Poissonization and other, new, technical devices.
 \\[.5cm]

\end{abstract}

\section{The structural distribution function}\label{strucdist}

Let the vector $X =(X_1,\dots,X_M)$ denote a mult$(n,p_M)$ distributed random
vector, where $p_M=(p_{M1} ,p_{M2}, \dots ,p_{MM})$ is the vector of
cell probabilities. Hence, the nonnegative components of $p_M$
satisfy $p_{M1} + \dots +p_{MM}=1$.

We will consider situations where $M=M_n$ is large with respect to
$n$, i.e.
\begin{equation}\label{nottozero}
M/n \not\to 0,\ \mbox{as}\ n\to\infty.
\end{equation}
 In these cases $X/n$
does not estimate $p_M$ accurately. For instance, for the average
mean squared error in estimating $Mp_{Mi}, i=1,\dots,M$, we have
$$ {1\over M}\sum_{i=1}^M \ex \Big(M\,  {X_i\over
n}-Mp_{Mi}\Big)^2
=
{M\over n}\sum_{i=1}^M p_{Mi}(1-p_{Mi})={M\over
n}\Big(1-\sum_{i=1}^n p_{Mi}^2\Big)\not\to 0, $$ unless
$\sum_{i=1}^M p_{Mi}^2 \to 1$ holds, i.e. unless $p_M$ comes close to a
unit vector $(0,\dots,0,1,0,\dots,0)$.

However, there are characteristics of $p_M$ that can be estimated
consistently. Here we will study
the {\em structural distribution function} of $p_M$. It is defined as
the empirical distribution function of the $Mp_{Mi},\ i=1,\dots M$,
and it  is given by
\begin{equation}
   F_M(x)=\frac{1}{M} \sum_{i=1}^{M} \I_{[Mp_{Mi}\leq x]},\ x\geq 0.
\end{equation}
Our basic assumption will be that $F_M$ converges weakly to a
limit distribution function $F$, i.e.
\begin{equation}\label{weakconv}
F_M \stackrel{w}{\to} F,\ \mbox{as}\ n\to\infty .
\end{equation}
The basic estimation problem is how to estimate $F_M$ (or $F$)
from an observation of $X$.

A  rule of thumb in statistics is to replace unknown probabilities
by sample fractions. This yields the so called {\em natural
estimator}. This estimator, denoted by  $\hat F_M$, is equal to
the empirical distribution function based on $M$ times the cell
fractions $X_i/n$, so
\begin{equation}
\hat F_M (x) = {1\over M} \sum_{i=1}^M \I_{[{M\over n}X_i\leq x]}.
\end{equation}
This estimator has often been used in linguistics, but turns out to be
inconsistent for estimating $F$; see Section \ref{inconsistency},
Khmaladze (1988), and Klaassen and Mnatsakanov (2000).

Our estimation  problem is related to   estimation in sparse multinomial
tables. For recent results on the estimation of cell probabilities in this
context see Aerts, Augustyns and Janssen (2000).

In Section \ref{simulation} we present a small simulation study of a typical
multinomial sample and the behavior of the natural estimator.
 It turns out that smoothing is required to obtain weakly consistent
estimators. An estimator based on grouping and an estimator based
on kernel smoothing are presented in Section \ref{smoothing}.
Section  \ref{technique} deals with the technique
of Poissonization and with the relation between weak and $L_1$ consistency.
These basic results are used in the weak
consistency proofs in Section \ref{consistency}.
Section \ref{discussion} contains a discussion.

\section{A simulation}\label{simulation}

\begin{figure}[h]
$$ \epsfxsize=6cm\epsfysize=4cm\epsfbox{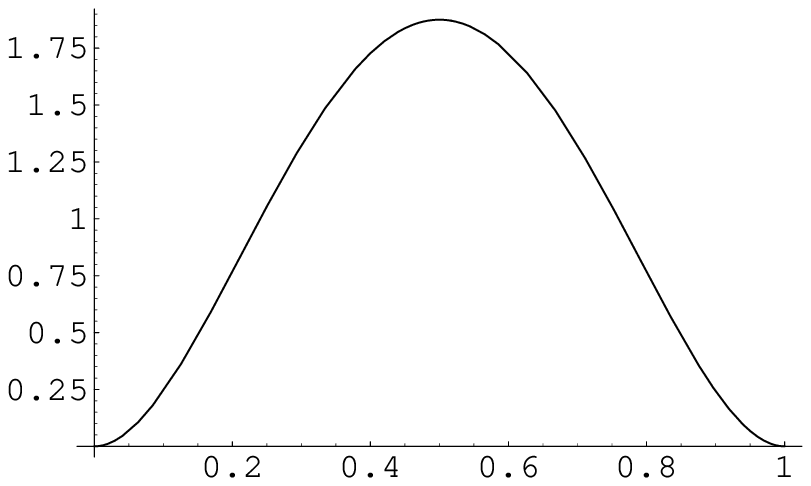}
 \epsfxsize=6cm\epsfysize=4cm\epsfbox{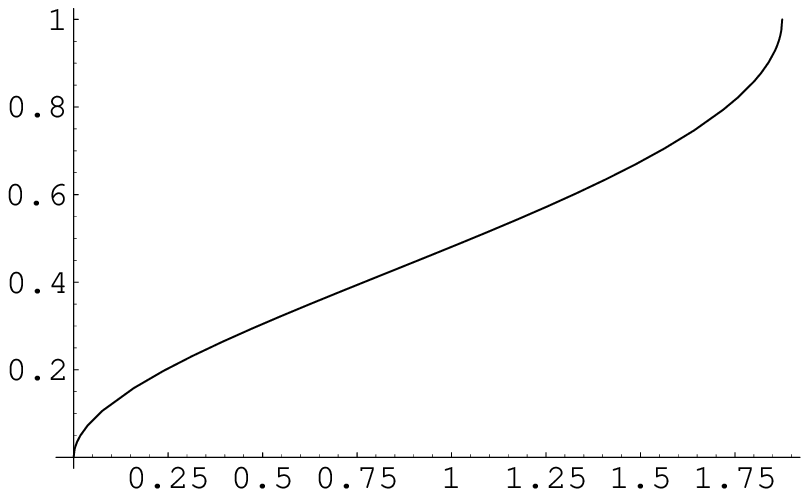} $$
\caption[]{The function $g$ and the corresponding structural distribution function $F$.
\label{fig:1}}
\end{figure}

We have simulated a sample with $M=1000$ and $n=2000$. The   cell
probabilities are generated via
\begin{equation}\label{probabilities}
p_{Mi}=G(i/M)-G((i-1)/M), i=1,\dots,M.
\end{equation}
 The
distribution function $G$ and its density $g$ have been chosen
equal to the functions
\begin{equation}
g(x)=30x^2(1 - x)^2\quad\mbox{and}\quad G(x)=10x^3 - 15 x^4 + 6 x^5,
0\leq x\leq 1.
\end{equation}
In Section \ref{smoothing} we show that   for these cell
probabilities, the limit structural distribution function $F$ from
(\ref{weakconv}) is equal to the distribution function of $g(U)$.
Here it is given by
\begin{equation}
F(x)=1-\sqrt{1-\sqrt{\frac{8}{15}x}},\quad 0\leq x \leq \frac{15}{8}.
\end{equation}
These functions are drawn in Figure \ref{fig:1}.

\begin{figure}[h]
$$ \epsfxsize=6cm\epsfysize=4cm\epsfbox{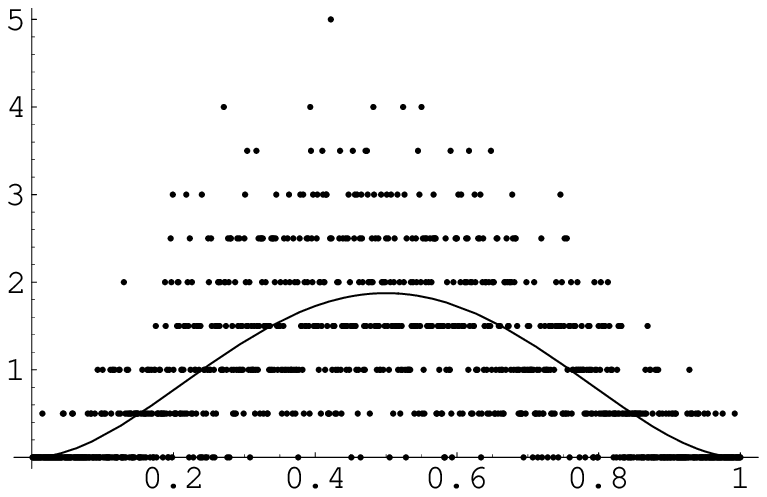}
 \epsfxsize=6cm\epsfysize=4cm\epsfbox{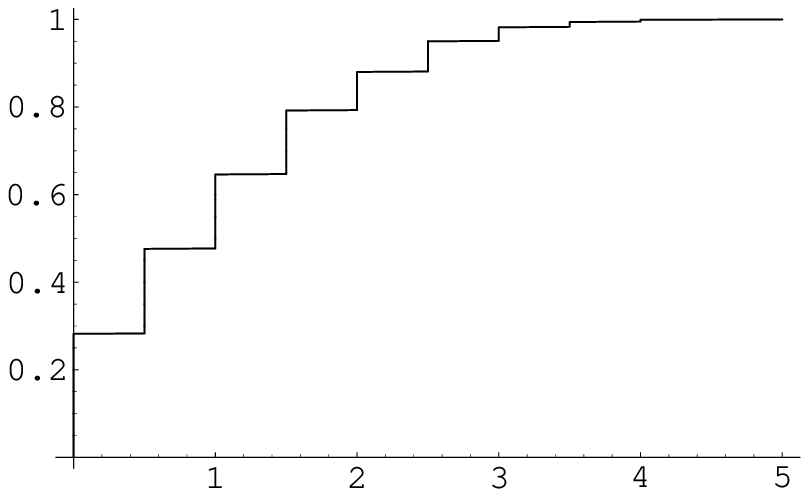} $$
\caption[]{The function $g$, $M/n$ times the cell counts, and in the second figure
 the natural estimator of $F$.
\label{fig:2}}
\end{figure}

\noindent For this simulated sample we have plotted the cell counts, multiplied by $M/n$,
and the natural estimate in Figure \ref{fig:2}. Comparison with the real $F$ in
Figure \ref{fig:1} clearly illustrates the inconsistency of the natural
estimator.

\section{Estimators based on smoothing techniques}\label{smoothing}

Up to now we have only assumed that the structural distribution function
$F_M$ converges weakly to a
limit distribution function $F$.
>From now on we will assume more structure.

Consider the function
\begin{equation}
g_M(u)=\sum_{i=1}^MMp_{Mi}\I_{({i-1\over M}, {i\over M}]}(u),\ u\in
 \mathbb R.
\end{equation}
This step function is a density
representing the cell probabilities and we shall call it the
 {\em parent density}.
The relation between this parent density $g_M$ and the structural
distribution function $F_M$ is given by the fact that if $U$ is a
uniform(0,1) random variable then $F_M$ is the distribution
function of $g_M(U)$. Note that
\begin{equation}
\ex g_M(U)=\int_{-\infty}^\infty g_M(u)du =\sum_{i=1}^M p_{Mi}= 1,
\end{equation}
so $g_M$ is a probability density indeed.

We will assume that there exists a limiting parent density $g$ on [0,1]
such that,
as $n\to\infty$,
\begin{equation}\label{condition}
\sup_{0<u\leq 1} |g_M(u)-g(u)|\to 0.
\end{equation}
Consequently we  have
$g_M(U)\to g(U)$, almost surely, and hence $F_M\stackrel{w}{\to} F$.

 The inconsistency of the natural estimator can be lifted by
first smoothing the cell counts $X_i$. We consider two smoothing
methods, grouping, which is actually some kind of histogram
smoothing, and a method based on  kernel smoothing of the counts.

\subsection{Grouping}

Let $m, k_j, j=0,1,\dots,m$, be integers, all depending on $n$, such that
$0=k_0<k_1<\dots < k_m=M$.
Define the group frequencies $\bar X_{j}$ as
\begin{equation}
\bar X_{j} = \sum_{i=k_{j-1}
+1}^{k_j}X_i,\quad j=1,\dots,m.
\end{equation}
Then the vector of grouped counts $\bar X$ is again multinomially distributed,
\begin{equation}
\bar X=(\bar X_1,\dots,\bar X_m)\sim mult (n,q_m),
\end{equation}
where $q_m=(q_{m1},\dots,q_{mm})$ and
\begin{equation}
q_{mj}=\sum_{i=k_{j-1}
+1}^{k_j}p_{Mi},\quad j=1,\dots,m.
\end{equation}
The grouped cells estimator, introduced in Klaassen and
Mnatsakanov (2000), is defined by
\begin{equation}
\hat F_M(x) = {1\over M} \sum_{j=1}^m(k_j-k_{j-1})
\I_{[{M\over{n(k_j-k_{j-1})}} \bar X_j\leq x]}, \quad x\geq 0.
\end{equation}
This estimator may be viewed as a structural distribution function with parent
density
\begin{equation}
\hat g_M(u)=\sum_{i=1}^m {M\over {n(k_i-k_{i-1})}}\,\bar X_i
\I_{[{k_{i-1}\over M}<u\leq{k_i\over M}]},\  u\in \mathbb R.
\end{equation}
This histogram is an estimator of the limiting parent density
 $g$ in (\ref{condition}). We
will prove weak consistency of the corresponding  estimator
$\hat F_M$ in Section
\ref{groupcons}.

\begin{figure}[h]
$$ \epsfxsize=6cm\epsfysize=4cm\epsfbox{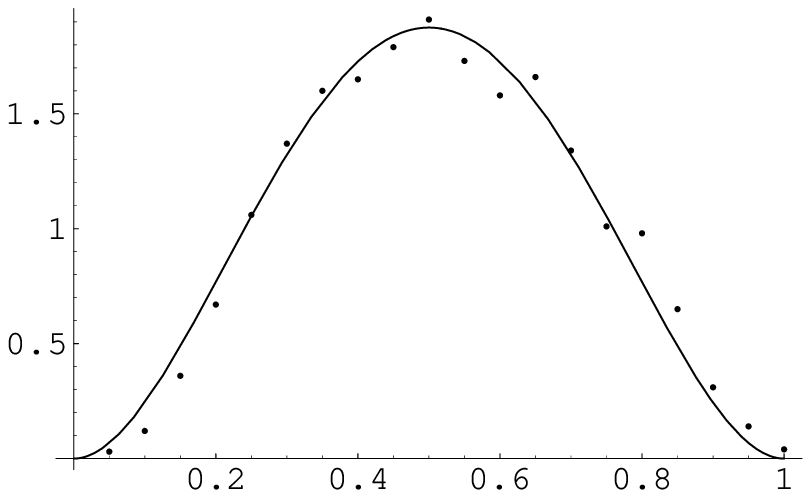}
 \epsfxsize=6cm\epsfysize=4cm\epsfbox{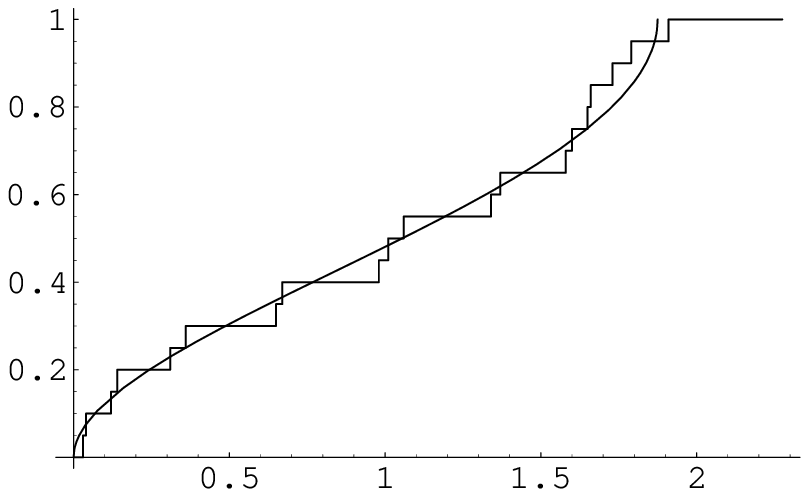} $$
\caption[]{$g, F$, and estimates   $\hat g_M$ and $\hat F_M$ by
 grouping with equal cell size.
\label{fig:3}}
\end{figure}

For our simulated example the estimates of $g$ and $F$ resulting
from grouping with equal group size $k=50$ are given in Figure
\ref{fig:3}.

\subsection{A kernel type estimator}

Now that we have seen that the estimator based on the grouped
cells counts is in fact
 based on a histogram estimate of the parent density $g$ we might also use
kernel smoothing to estimate $g$ and proceed in a similar manner.
If we choose a  probability density $w$ as {\em kernel function} and
a {\em bandwidth} $k\geq 0$, we get the following estimator for the parent
 density $g$
\begin{equation}
\hat g_M(u)={M\over nk}\sum_{i=1}^M w\Big({\lceil Mu
\rceil-i\over k}\Big) X_i
,\ u\in \mathbb R.
\end{equation}
As an estimator for the structural distribution function of the function
$F$ we take the empirical distribution function of $\hat g_M(U)$ with $U$
uniform, namely
\begin{equation}
\hat F_M(x)={1\over M}\sum_{j=1}^M \I_{[{M\over nk}\sum_{i=1}^M
w({j-i\over k} )X_i\leq
x]}.
\end{equation}
Weak consistency of this estimator will be derived in Section
\ref{kerncons}.

\begin{figure}[h]
$$ \epsfxsize=6cm\epsfysize=4cm\epsfbox{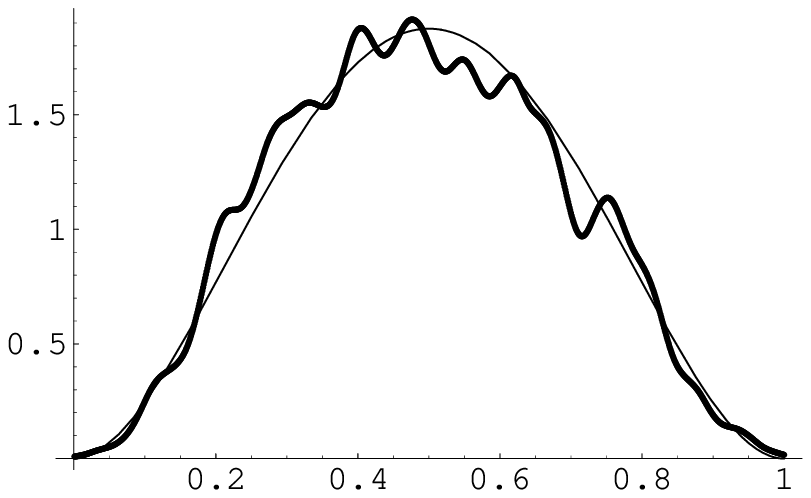}
 \epsfxsize=6cm\epsfysize=4cm\epsfbox{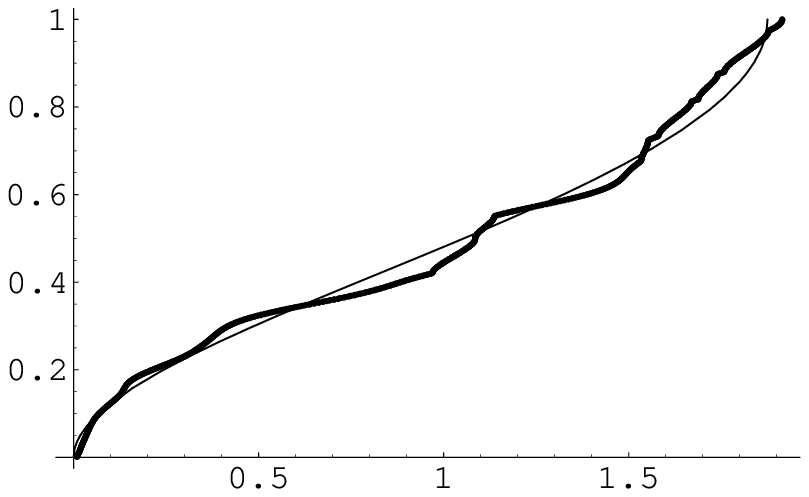} $$
\caption[]{$g, F$, and estimates  $\hat g_M$ and $\hat F_M$
 by kernel smoothing.
\label{fig:4}}
\end{figure}

For our simulated example kernel estimates $\hat g_M$ and $\hat F_M$
of $g$ and $F$, respectively,
with $k$ equal to 50 are given in Figure \ref{fig:4}.

\section{Relevant techniques}\label{technique}

In our proofs we shall use repeatedly the powerful method of
Poissonization and a device involving $L_1$ convergence.

\subsection{Poissonization}\label{Poissonization}

Consider the random vectors $X$ and $Y$, with
\begin{equation}\label{pois}
X=(X_1,\dots,X_M)\sim \mbox{mult}(n,p_M)
\quad\mbox{and}\quad
Y=(Y_1,\dots,Y_M),\ Y_i\sim \mbox{Poisson}(np_{Mi}),
\end{equation}
where $Y_1,\dots,Y_M$ are independent.
Note
\begin{equation}
N=\sum_{i=1}^M Y_i \sim \mbox{Poisson}(n).
\end{equation}
Given $N=k$ the random vector $Y$ has a mult$(k,p_M)$ distribution.

Based on an infinite sequence of $mult(1,p_{M1},\dots,p_{MM})$ random vectors
one can construct vectors $X$ and $Y$, the cell counts over $n$ and $N$
of these vectors repectively, with the distributions (\ref{pois}).
Given $N=k$ they are coupled as follows

\parbox{11cm}{\begin{eqnarray*}
k\leq n: & X=Y + \mbox{mult}(n-k,p_M),\\ k>n: &
Y=X+\mbox{mult}(k-n,p_M).
\end{eqnarray*}
}
\hfill
\parbox{2cm}{\begin{equation}\label{condpois}\end{equation}}

\noindent Note that this shows that either $X_i\leq Y_i$ for all $i$ or
$X_i\geq Y_i$ for
all $i$.

\subsection{Convergence in $L_1$ and weak convergence}\label{L1}

An important step in the (in)consistency proofs is to show that
``Poissonization
is allowed'', i.e. that we can transfer the limit result for the estimator
based on the Poissonized sample, the ``Poissonized version'',
 to the original estimator. The following proposition is used repeatedly,
also if no Poissonized version is involved.

\begin{prop}\label{prop:1}
Let $F$ be a distribution function and let $\hat F_n$ and $\tilde
F_n$ be possibly random distribution functions.
If
\begin{equation}
\label{ass:1}
\tilde F_n \stackrel{w}{\to} F,\quad\mbox{in probability},
\end{equation}
and
\begin{equation}
\label{ass:2}
\int |\hat F_n - \tilde F_n|\stackrel{P}{\to}0
\end{equation}
hold, then
\begin{equation}
\hat F_n \stackrel{w}{\to} F,\quad\mbox{in probability},
\end{equation}
is valid, i.e. for all $\epsilon >0$ and all continuity points $x_0$ of $F$
\begin{equation}
P(|\hat F_n(x_0) -   F(x_0)|>\epsilon)\to 0.
\end{equation}
\end{prop}
In  the special case where $ \tilde F_n$ equals $F$,
the proposition states that $L_1$ convergence implies weak
convergence.

\bigskip

\noindent{\bf Proof}
Note that for all $x_0$ and all $\delta>0$ we have
\begin{equation}\label{int}
\int_{x_0-\delta}^{x_0+\delta}|\hat F_n -  F | \leq
\int_{-\infty}^{\infty}|\hat F_n -  \tilde F |  +
\int_{x_0-\delta}^{x_0+\delta}|\tilde F_n -  F |.
\end{equation}

 Let $x_0$ denote an arbitrary
continuity point of $F$ and $\epsilon$ an arbitrary positive number.
Choose $\delta>0$ such that $F(x_0+\delta)-F(x_0-\delta)\leq \epsilon$
and such that
$x_0-\delta $ and  $x_0+\delta$
are continuity points of $F$.
Then
\begin{equation}
|\tilde F_n(x_0-\delta)-F(x_0-\delta)|<\epsilon\ \mbox{and}\
|\tilde F_n(x_0+\delta)-F(x_0+\delta)|<\epsilon
\end{equation}
imply
\begin{equation}
\int_{x_0-\delta}^{x_0+\delta}|\tilde F_n -  F |<4\delta\epsilon.
\end{equation}
Hence, we have

\parbox{11cm}{\begin{eqnarray*}
\lefteqn{P\Big(\int_{x_0-\delta}^{x_0+\delta}|\tilde F_n -  F |
\geq 4\delta\epsilon\Big)}\\ &\leq& P(|\tilde
F_n(x_0-\delta)-F(x_0-\delta)|\geq\epsilon) + P(|\tilde
F_n(x_0+\delta)-F(x_0+\delta)|\geq\epsilon) \nonumber
\end{eqnarray*}}\hfill\parbox{2cm}{\begin{equation}\end{equation}}

\noindent and, by (\ref{ass:1}),
\begin{equation}
\int_{x_0-\delta}^{x_0+\delta}|\tilde F_n -  F |\stackrel{P}{\to} 0.
\end{equation}
Consequently, by (\ref{ass:2}) and (\ref{int}) we get
\begin{equation}
\int_{x_0-\delta}^{x_0+\delta}|\hat F_n -  F |\stackrel{P}{\to} 0.
\end{equation}

Choose $0<\delta'<\delta$ such that $F(x_0+\delta')\leq F(x_0)+\frac{1}{2}
\epsilon$
and $F(x_0-\delta')\geq F(x_0)-\frac{1}{2} \epsilon$. Then we see
\begin{equation}
|\hat F_n(x_0)-F(x_0)|\geq\epsilon \Rightarrow
\int_{x_0-\delta'}^{x_0+\delta'}|\hat F_n -  F | \geq
\frac{1}{2}\delta' \epsilon
\end{equation}
and hence
\begin{eqnarray*}
\lefteqn{ P(|\hat F_n(x_0)-F(x_0)|\geq\epsilon)\leq
P\Big(\int_{x_0-\delta'}^{x_0+\delta'}|\hat F_n -  F | \geq
 \frac{1}{2}\delta' \epsilon\Big)}\\
&\leq&
P\Big(\int_{x_0-\delta}^{x_0+\delta}|\hat F_n -  F | \geq
 \frac{1}{2}\delta' \epsilon\Big)\to 0.
\end{eqnarray*}
Since this holds for arbitrary continuity points $x_0$ and
arbitrary $\epsilon>0$ we have established $\hat
F_n\stackrel{w}{\to}F$, in probability. \hfill$\Box$

\section{Consistency}\label{consistency}

\subsection{The natural estimator}\label{inconsistency}

The basic trick in dealing with the difference of the natural estimator
and its Poissonized version,
\begin{equation}
\tilde F_M(x)={1\over M} \sum_{i=1}^M \I_{[M {Y_i\over n}\leq x]},
\end{equation}
uses the  coupling   as in (\ref{condpois}) and is given by
 the following string of inequalities

\parbox{11cm}{\begin{eqnarray*}
\lefteqn{|\hat F_M(x) - \tilde F_M(x)|
\leq {1\over M} \sum_{i=1}^M |\I_{[M {X_i\over n}\leq x]}-
\I_{[M {Y_i\over n}\leq x]}|}\\
&\leq {1\over M} \sum_{i=1}^M \I_{[ X_i\not= Y_i]}
\leq {|N-n|\over M}=O_P\Big({\sqrt{n}\over M}\Big).
\end{eqnarray*}}\hfill\parbox{2cm}{\begin{equation}\label{poisallowed}
\end{equation}}

\noindent By (\ref{nottozero}) the right hand side converges to zero in
probability
and this shows that Poissonization is allowed.

Because of the independence of the Poisson counts $Y_i$
we can easily bound the variance of the Poissonized estimator. We get
\begin{equation}
\var \tilde F_M(x) = \var {1\over M} \sum_{i=1}^M \I_{[M {Y_i\over n}\leq x]}
\leq {1\over 4M} \to 0.
\end{equation}
We also have
\begin{equation}
\ex \tilde F_M(x) = {1\over M}\sum_{i=1}^M P({M\over n} Y_i\leq x)
\not= {1\over M}\sum_{i=1}^M   \I_{[M p_{Mi}\leq x]}=F_M(x)
\end{equation}
and

\parbox{13cm}{\begin{eqnarray*}
\lefteqn{\ex \int x^2d\tilde F_M(x)=\ex {1\over M} \sum_{i=1}^M
\Big({M\over n} Y_i\Big)^2}\\
&=
{1\over M} \sum_{i=1}^M \Big({M\over n}\Big)^2\{np_{Mi} + (np_{Mi})^2\}
={M\over n} + \int  x^2d\  F_M(x).
\end{eqnarray*}}\hfill\parbox{2cm}{\begin{equation}\end{equation}}\\
Together with (\ref{nottozero})
this gives two reasons why $\tilde F_M(x)$ is probably not a consistent
estimator of $F$.
Then, by (\ref{poisallowed}) the natural estimator has to
be inconsistent too.

The inconsistency of the
structural distribution function  has been
established in Khmaladze (1988), Khmaladze and Chitashvili (1989),
Klaassen and Mnatsakanov (2000) and Van Es and Kolios (2002).
In these papers the situation is considered of
a {\em large number of rare events}, i.e.  $n/M\to\lambda$ for some constant
$\lambda$.
The explicit limit in probability of $\hat F_M(x)$ turns out to be a Poisson mixture
 of $F$ then.

\subsection{Grouping}\label{groupcons}

Under the additional assumption $n/M\to\lambda$, for some constant $\lambda$,
weak consistency of the estimator based on grouped cells has been
proved, without using Poissonization, by Klaassen and Mnatsakanov (2000) and by the
Poissonization method for the simpler case of equal group size,
i.e. $k_j=k$, by Van Es and Kolios (2002).
We shall prove the following generalization without using Poissonization.

\begin{thm}\label{thm:group}
If $m/n\to 0$,
\begin{equation}
\sup_{1\leq j\leq m} {{k_j-k_{j-1}}\over M}\to 0,
\end{equation}
and
\begin{equation}
\sup_{0<u\leq 1} |g_M(u)-g(u)|\to 0
\end{equation}
are valid for some limiting parent density
$g$ that  is continuous on $[0,1]$, then
\begin{equation}
\hat F_M \stackrel{w}{\to}F,\quad\mbox{in probability},
\end{equation}
holds with
\begin{equation}
\hat F_M(x)={1\over M}\sum_{j=1}^m (k_j-k_{j-1}) \I_{[{M\over
n(k_j-k_{j-1})}\sum_{i=k_{j-1}+1}^{k_{j}} X_i\leq x]}.
\end{equation}
\end{thm}

\noindent {\bf Proof}

The estimator $\hat F_M$ behaves asymptotically as
\begin{equation}
\bar F_M(x) = {1\over M}\sum_{j=1}^m(k_j-k_{j-1})\I_{[{Mq_{mj}\over k_j-k_{j-1}}\leq x
]}.
\end{equation}
Indeed, in view of $\int |\I_{[a\leq x]}-\I_{[b\leq
x]}|dx=|b-a|$ we have
\begin{eqnarray}
\lefteqn{\int |\hat F_M(x) - \bar F_M(x)|dx } \nonumber\\ &\leq&
\int \sum_{j=1}^m {k_j-k_{j-1}\over M} \Big|\I_{[{M\bar X_j\over
n(k_j-k_{j-1})} \leq x ]} - \I_{[{Mq_{mj}\over k_j-k_{j-1}}\leq x
]}  \Big|dx\\ &=& \sum_{j=1}^m {k_j-k_{j-1}\over M} \Big|{M\bar
X_j\over n(k_j-k_{j-1})}
 - {Mq_{mj}\over k_j-k_{j-1}} \Big|. \nonumber
\end{eqnarray}
Consequently, we obtain
\begin{eqnarray*}
\lefteqn{\ex \int |\hat F_M(x) - \bar F_M(x)|dx
\leq {m\over n}\, {1\over m} \sum_{j=1}^m \ex |  \bar X_j - n q_{mj}|}\\
&\leq&
{m\over n}\sqrt{ {1\over m} \sum_{j=1}^m \ex (  \bar X_j - n
q_{mj})^2}
=
{m\over n}\sqrt{ {1\over m} \sum_{j=1}^m nq_{mj}(1-q_{mj})}\\
&\leq&
 \sqrt{{m\over n}}\to 0
\end{eqnarray*}
and hence
\begin{equation}
 \int |\hat F_M(x) - \bar F_M(x)|dx \stackrel{P}{\to} 0.
\end{equation}
In order to prove $\hat F_M\stackrel{w}{\to} F$ in probability, by
Proposition \ref{prop:1} it remains to show $\bar
F_M\stackrel{w}{\to} F$.

Consider the function
\begin{equation}
\bar g_M(u)=\sum_{j=1}^m {1\over
{k_j-k_{j-1}}}\sum_{i=k_{j-1}+1}^{k_j}
Mp_{Mi}\I_{({k_{j-1}\over M},{k_j\over M}]}(u).
\end{equation}
For ${k_{j-1}/ M}<u\leq{k_j/M}$ we have
\begin{eqnarray*}
\lefteqn{|\bar g_M(u) -g(u)|\leq {1\over
{k_j-k_{j-1}}}\sum_{i=k_{j-1}+1}^{k_j} |Mp_{Mi}-g(u)|}\\
&\leq&
\sup_{{k_{j-1}/ M}<v\leq{k_j/ M}} |g_M(v) - g(u)|\\
&\leq&
\sup_v |g_M(v) - g(v)|
+ \sup_{|u-v|\leq \sup_j (k_j - k_{j-1})/M}|g(v)-g(u)|.\\
\end{eqnarray*}
By assumption, the function
$g$ is uniformly continuous and hence $\sup_j (k_j
- k_{j-1})/M\to 0$ implies $\bar g_M(U)\to g(U)$, almost surely,
and in distribution, i.e. $\bar F\stackrel{w}{\to} F$, which
completes the proof of the theorem. \hfill$\Box$

\subsection{The kernel type estimator}\label{kerncons}

Weak consistency of the kernel type estimator is established by
the next theorem.

\begin{thm}\label{thm:kern}
If $k\to \infty, {k/M}\to 0, {M/(nk)}\to 0$ hold,   if
$w$ is a density that is Riemann integrable on bounded intervals,
that is also Riemann square
integrable on bounded intervals, and that
has   bounded support or is ultimately monotone in its tails, and if
\begin{equation}\label{w4}
\sup_{0<u\leq 1} |g_M(u)-g(u)|\to 0
\end{equation}
holds with $g$   continuous on $[0,1]$, then
\begin{equation}
\hat F_M \stackrel{w}{\to}F,\quad\mbox{in probability},
\end{equation}
is valid for
\begin{equation}
\hat F_M(x)={1\over M}\sum_{j=1}^M \I_{[{M\over nk}\sum_{i=1}^M
 w({j-i\over k})X_i\leq
x]}.
\end{equation}
\end{thm}

\noindent{\bf Proof}
Let
\begin{equation}
\tilde F_M(x)={1\over M}\sum_{j=1}^M \I_{[{M\over nk}\sum_{i=1}^M
 w({j-i\over k})Y_i\leq
x]}
\end{equation}
be the Poissonized version of $\hat F_M(x)$.
Note that by the coupling argument $X_i\geq Y_i$ for all $i$ or
$X_i\leq Y_i$ for all $i$. Since $w$ is a Riemann integrable density
we thus get
\begin{eqnarray*}
\lefteqn{\ex \int |\hat F_M(x) - \tilde F_M(x)|dx}\\
&\leq&
\ex {1\over M} \sum_{j=1}^M\int |\I_{[{M\over nk}\sum_{i=1}^M
 w({j-i\over k})X_i\leq
x]}-\I_{[{M\over nk}\sum_{i=1}^M
 w({j-i\over k})Y_i\leq x]}|dx\\
&=&
\ex {1\over M} \sum_{j=1}^M\Big|{M\over nk}\sum_{i=1}^M
 w\Big({j-i\over k}\Big)(X_i-Y_i)\Big|
=
\ex {1\over M} \sum_{j=1}^M{M\over nk}\sum_{i=1}^M
 w\Big({j-i\over k}\Big)|X_i-Y_i |\\
&=&
\ex{1\over n}\sum_{i=1}^M\Big( \sum_{j=1}^M{1\over k}\,w\Big({j-i\over
k}\Big)\Big)|X_i-Y_i |
\leq
 \sum_{l\in {\mathbb Z}}{1\over k}\,w\Big({l\over k}\Big)
\ \ex {|N-n|\over n} = O\Big({1\over \sqrt{n}}\Big).
\end{eqnarray*}
Consequently, by Proposition \ref{prop:1} it suffices to prove
\begin{equation}\label{w}
\tilde F_M\stackrel{w}{\to} F, \quad\mbox{in probability}.
\end{equation}
Define
\begin{equation}
\bar F_M(x)={1\over M}\sum_{j=1}^M \I_{[{1\over k}\sum_{i=1}^M
 w({j-i\over k})Mp_{Mi}\leq
x]}.
\end{equation}
To prove (\ref{w}), by Proposition \ref{prop:1}, it suffices to prove
\begin{equation}\label{w2}
\ex \int |\tilde F_M(x) - \bar F_M(x)|dx\stackrel{P}{\to} 0\quad\mbox{and}
\quad \bar F_M\stackrel{w}{\to} F, \ \mbox{in probability}.
\end{equation}
Indeed, since the $Y_i$ are independent and $w$ is square Riemann integrable, we have
\begin{eqnarray*}
\lefteqn{\ex \int |\tilde F_M(x) - \bar F_M(x)|dx
\leq {1\over M} \sum_{j=1}^M\ex\Big|
{M\over nk}\sum_{i=1}^M w\Big({j-i\over
k}\Big)(Y_i-np_{Mi})\Big|}\\
&\leq&
\sqrt{{1\over M} \sum_{j=1}^M\var\Big\{
{M\over nk}\sum_{i=1}^M w\Big({j-i\over
k}\Big)(Y_i-np_{Mi})\Big\}}\\
&=&
\sqrt{{M\over n^2k^2} \sum_{j=1}^M
\sum_{i=1}^M w^2\Big({j-i\over
k}\Big) np_{Mi}}
\leq
\sqrt{{M\over nk^2} \sum_{i=1}^M
\sum_{\ell\in \mathbb Z} w^2\Big({\ell\over
k}\Big) p_{Mi}}\\
&=&
\sqrt{{M\over nk}
\sum_{\ell\in \mathbb Z} {1\over k}\,w^2\Big({\ell\over
k}\Big)}
=
O\Big(\sqrt{{M\over nk}}\Big)\to 0,
\end{eqnarray*}
because of $k\to \infty$ and $M/(nk)\to 0$.
This proves the first statement of (\ref{w2}).

Finally, we prove the second statement of (\ref{w2}).
As parent density for the distribution function $\bar F_M$ we choose
\begin{equation}
\bar g_M(u)=\sum_{j=1}^M\ {1\over k}\sum_{i=1}^M
w\Big({j-i\over k}\Big)Mp_{Mi}\ \I_{({{j-1}\over M},{j\over M}]}(u),
\ u \in \mathbb R.
\end{equation}

Note that $g_M$ vanishes outside (0,1].
Fix $u\in(0,1)$. For $u\in({{j-1}\over M},{j\over M}]$, and $K>0$ fixed,
 we have
\begin{eqnarray*}
\lefteqn{
|\bar g_M(u)-  g(u)|}\nonumber\\
&\leq&
\sum_{\ell\in \mathbb Z} {1\over k}\,w\Big({\ell\over k}\Big)
\Big|g_M\Big({j-\ell\over M} \Big)
-g\Big({j-\ell\over M}\Big)\Big|\nonumber\\
&+&
\sum_{|\ell|\leq Kk} {1\over k}\,w\Big({\ell\over k}\Big)
\Big|g \Big({j-\ell\over M}\Big)
-g(u)\Big|\\
&+&
\Big\{
\sum_{|\ell|> Kk} {1\over k}\,w\Big({\ell\over k}\Big)+
\Big|\sum_{\ell\in \mathbb Z} {1\over k}\,w\Big({\ell\over
k}\Big)-1\Big|\Big\}\sup_u g(u).
\nonumber
\end{eqnarray*}
Note that the conditions imposed on $w$ guarantee
that
\begin{equation}
\sum_{|\ell|> Kk} {1\over k}\,w\Big({\ell\over k}\Big)
\end{equation}
 is arbitrarily
small for $K$ sufficiently large,
 that
\begin{equation}
\sum_{|\ell|\leq Kk} {1\over k}\,w\Big({\ell\over k}\Big)
\to \int_{-K}^K w(u)du,
\end{equation}
which is arbitrarily close to one for $K$ large enough, and hence that
\begin{equation}
\sum_{\ell\in \mathbb Z} {1\over k}\,w\Big({\ell\over
k}\Big)\to 1,
\end{equation}
as $k\to \infty$. Consequently, in view of
(\ref{w4}), and in view of the uniform continuity and boundedness of $g$,
 all three terms at the right hand side tend to zero as $k\to\infty$ and
subsequently $K\to\infty$.
So, $\bar g_M(U)\to g(U)$, almost surely and in distribution, which
implies $\bar F_M\stackrel{w}{\to} F$.

\hfill$\Box$

\section{Discussion}\label{discussion}

The key assumption in the consistency proofs of the
grouping and kernel estimators is the existence of
a limiting parent density. This is a reasonable assumption
only, if there is a natural ordering of the cells and
neighboring cells have approximately the same
cell probabilities. In applications like e.g. linguistics
this need not be the case. Consider a text of $n$ words
of an author with a vocabulary of $M$ words. Here the words
in the vocabulary correspond to the cells of the multinomial
distribution and the existence of a limiting or approximating
parent density is rather unrealistic. To a lesser extent
this might be the case in biology, where cells correspond to
species and $n$ is the number of individuals found in some
ecological entity.

An estimator that is consistent even if our key assumption
does not hold, has been constructed in Klaassen and
Mnatsakanov (2000). However, it seems to have a logarithmic
rate of convergence only. The rates of convergence of our
grouping and kernel estimators will depend on the rate
at which the assumed limiting parent density can be estimated.
This issue is still to be investigated, but under the assumption
$n/M\to\lambda$, for some constant
$\lambda$,
 Van Es and Kolios (2002) show that,
for the relatively simple case of equal group size, an algebraic
rate of convergence can be achieved by the estimator based on
grouping.

Since the estimators studied here are based on smoothing of the cell
frequencies an important open problem is the choice of the
smoothing parameter. For the estimator based on grouping this is
the choice of  the sizes of the groups and for the kernel type
estimator the choice of the bandwidth.
By studying convergence rates these choices may be optimized.

\vspace{2cm}

\noindent{\bf\Large Acknowledgement}

\medskip

This paper has been prepared under
INTAS-97-Georgia-1828.

\bigskip

\noindent{\bf\Large References}

\begin{description}

\item

Aerts, M., Augustyns, I. and P. Janssen (2000), Central limit
theorem for the total squared error of local polynimial estimators
of cell probabilities, {\em J. Statist. Plann. Inference}, 91,
181--193.

\item
Van Es, B. and S. Kolios (2002), { Estimating a structural distribution
 function by grouping}, Mathematics ArXiv PR/0203080.

\item
Khmaladze, E.V. (1988), The statistical analysis of a large number of rare
events, Report MS-R8804, CWI, Amsterdam.

\item
Khmaladze, E.V. and R.Ya. Chitashvili (1989),  Statistical analysis of
a large number of rare events
and related problems (Russian),
Proc. A. Razmadze Math. Inst. Georgian Acad. Sci., Tbilisi, 92, 196-245.

\item
Klaassen, C.A.J. and R.M. Mnatsakanov (2000),   Consistent
estimation of the structural distribution function, {\em Scand.
J. Statist.}, 27, 733--746.

\end{description}

\end{document}